\documentclass[11pt]{article}

\usepackage{amsmath}
\usepackage{amssymb}
\usepackage{eucal}

\begin{document}

\baselineskip 16pt

\title{On a lattice characterization of finite soluble $PST$-groups}

\author{Zhang Chi \thanks{Research of the first author is supported by
 China Scholarship Council and NNSF of
China(11771409)}\\
{\small Department of Mathematics, University of Science and
Technology of China,}\\ {\small Hefei 230026, P. R. China}\\
{\small E-mail:
zcqxj32@mail.ustc.edu.cn}\\ \\
{ Alexander  N. Skiba}\\
{\small Department of Mathematics and Technologies of Programming,
  Francisk Skorina Gomel State University,}\\
{\small Gomel 246019, Belarus}\\
{\small E-mail: alexander.skiba49@gmail.com}}

\date{}
\maketitle

\begin{abstract} Let   $\mathfrak{F}$ be  a class of finite  groups   and  $G$   a
 finite group. Let  ${\cal L}_{\mathfrak{F}}(G)$ be  the
set of all subgroups $A$ of $G$ with $A^{G}/A_{G}\in \mathfrak{F}$.
  A chief factor $H/K$ of $G$
 is \emph{ $\mathfrak{F}$-central} in
$G$  if  $(H/K)\rtimes (G/C_{G}(H/K)) \in\mathfrak{F}$.

We  study  the structure of $G$ under the hypothesis that      every chief factor
 of $G$ between $A_{G}$ and $A^{G}$ is
$\mathfrak{F}$-central in $G$ for every subgroup $A\in {\cal 
L}_{\mathfrak{F}}(G)$.  As an application, we prove that
a finite soluble group $G$  is a $PST$-group if and only if
$A^{G}/A_{G}\leq Z_{\infty}(G/A_{G})$ for every subgroup $A\in {\cal 
L}_{\mathfrak{N}}(G)$, where  $\mathfrak{N}$ is the class of all
nilpotent groups.

\end{abstract}

\let\thefootnoteorig\thefootnote
\renewcommand{\thefootnote}{\empty}

\footnotetext{Keywords: finite group,  lattice of subgroups,
 $\mathfrak{F}$-central
 chief factor,  saturated formation, Fitting formation.}

\footnotetext{Mathematics Subject Classification (2010): 20D10,
20D15, 20D20, 20E15}
\let\thefootnote\thefootnoteorig

\section{Introduction}

Throughout this paper, all groups are finite and $G$ always denotes
a finite group. Moreover,  ${\cal L}(G)$  denotes the lattice of all
subgroups of $G$.  We use $A^{G}$ to denote the normal closure of the subgroup $A$ in $G$; 
$A_{G}=\bigcap _{x\in G} A^{x}$.  If $L\leq T$ are normal subgroups of $G$, then we say that $T/L$ is a 
\emph{normal section} of $G$.
 In this paper,
 $\mathfrak {F}$ is  a class of groups containing all identity groups;
 $\mathfrak{N}$  denotes the classes of all
nilpotent  groups.

    Wielandt proved \cite{WI} that  the set ${\cal L}_{sn}(G)$, of
 all subnormal subgroups of a finite group $G$,  forms a 
sublattice   of the lattice  ${\cal L}(G)$.  Later, Kegel proposed in 
\cite{5}   an original  idea of generalization the lattice ${\cal L}_{sn}(G)$
  being based of the  group classes theory.
The papers \cite{WI, 5} has become a motivation for a lot of others  studies related to 
finding and applying sublattices of the lattices ${\cal L}(G)$ and ${\cal L}_{sn}(G)$
  (see, for example,
 \cite{Keg, 19, 20},  Ch.
 6 in \cite{15} and the recent paper \cite{1}).

In this paper, we discuss a new approach that   allows to allocate two new
 classes of sublattices   in the lattice ${\cal L}(G)$ and give some  applications
 of such   sublattices in the theory of generalized $T$-groups.

Let   $\Delta$ be  any set of normal sections of $G$. Then we say 
 that $\Delta$ is \emph{$G$-closed} provided for any two 
$G$-isomorphic normal sections  $H/K$ and $T/L$, where $T/L\in  \Delta$, we have  
 $H/K\in  \Delta$. If $L\leq T$ are normal subgroups of $G$,
 then we write  
$T/L\leq Z_{\Delta}(G)$ (or simply $T\leq Z_{\Delta} (G)$ if $L=1$) 
provided either $L=T$ or $L < T$ and $H/K\in \Delta$ for 
 every  chief factor $H/K$ of $G$ between $L$ and $T$.

 Before  continuing, we recall some notations and concepts of the  group classes  theory.  
The symbol  
$G^{\mathfrak {F}}$ denotes the \emph{$\mathfrak {F}$-residual} of $G$, that is,
 the intersection of all normal subgroups 
$N$ of $G$ with  $G/N \in \mathfrak {F}$; $G_
{\mathfrak {F}}$ denotes the \emph{$\mathfrak {F}$-radical} of $G$, that is,
 the   product of all normal subgroups 
$N$ of $G$ with  $N \in \mathfrak {F}$. The class $\mathfrak {F}$ 
is said 
to be: \emph{normally  hereditary}   if $H\in {\frak{F}}$
 whenever   $ H \trianglelefteq  G 
\in {\frak{F}}$;  \emph{saturated}  if $G\in {\frak{F}}$
whenever $G^{\frak{F}}\leq \Phi (G)$; 
  a {\sl formation} if  every homomorphic image of  $G/G^{\mathfrak {F}}$
belongs to $\mathfrak {F}$ for any group $G$; a {\sl Fitting class } if  every normal subgroup 
 of  $G_{\mathfrak {F}}$
belongs to $\mathfrak {F}$ for any group $G$.

Now, let    ${\cal L}_{\Delta}(G)$ be  the set of all
subgroups $A$ of $G$ such that  $A^{G}/A_{G}\leq Z_{\Delta}(G)$, and  
 we put ${\cal L}_{\mathfrak{F}}(G)={\cal L}_{\Sigma }(G)$ if $\Sigma $ is 
the set of all normal sections $T/L$  of $G$ with  $T/L\in \mathfrak{F}$.

{\bf Theorem  1.1.}  (i) {\sl The set ${\cal L}_{\Delta}(G)$ is a
  sublattice of the lattice ${\cal L}(G)$ for any $G$-closed set of chief factors of $G$.  }

(ii) {\sl If $\mathfrak{F}$ is a normally hereditary formation, then the
 set ${\cal L}_{\mathfrak{F}}(G)$ is a lattice (a meet-sublattice of ${\cal L}(G)$} 
 \cite[p. 7]{Schm}).

(iii) {\sl If $\mathfrak{F}$ is a Fitting formation, then
 ${\cal L}_{\mathfrak{F}}(G)$ is a  sublattice  of the lattice ${\cal L}(G)$.}

A  subgroup $M$ of  $G$ is called  \emph{ modular} in $G$ if $M$ is a modular element
(in the sense of  Kurosh \cite[ p. 43]{Schm})  of the
  lattice ${\cal L}(G) $. It is known that \cite [Theorem 5.2.3]{Schm} for 
every  modular subgroup $A$ of $G$ all chief factors of $G$ between 
$A_{G}$ and $A^{G}$ are cyclic.
Consequently, despite the fact that in the general case the intersection 
of two modular subgroups of $G$ may  be non-modular, the following result 
holds.

{\bf Corollary    1.2.}    {\sl If $A$ and $ B$ are modular subgroups of $G$, then
 every chief factor of $G$ between $(A\cap B)_{G}$     and $(A\cap B)^{G}$ is cyclic.      }

   A subgroup $A$ of $G$ is said to be 
\emph{quasinormal }  (respectively \emph{$S$-quasinormal }  or
\emph{$S$-permutable} \cite{prod}) in $G$  if $A$
 permutes with all subgroups (respectively
  with all Sylow subgroups) $H$
of $G$, that is, $AH=HA$.  For every quasinormal subgroup $A$ of $G$ we 
have  $A^{G}/A_{G} \leq Z_{\infty}(G/A_{G})$ \cite[Corollary 1.5.6]{prod}. It is known also that, in general,  
the intersection  quasinormal   subgroups   of $G$ may be 
 non-quasinormal. Nevertheless, the following fact holds.

{\bf Corollary    1.3.}    {\sl If $A$ and $ B$ are quasinormal  subgroups of $G$, then
 $(A\cap B)^{G}/(A\cap B)_{G} \leq Z_{\infty}(G/(A\cap B)_{G})$. }

Recall that a chief factor $H/K$ of $G$ is said to be 
  \emph{ $\mathfrak{F}$-central} in
$G$   \cite{14}  if  $(H/K)\rtimes (G/C_{G}(H/K)) \in \mathfrak{F} $.

Our next result is the following 

{\bf Theorem  1.4.} {\sl  Let  $\mathfrak{F}$ be  a normally hereditary
 saturated formation containing all nilpotent groups and  
 $\Delta$ the set of all 
$\mathfrak{F}$-central chef factors of $G$. }

(i) {\sl   If 
the $\mathfrak{F}$-residual $D=G^{\mathfrak{F}}$ of $G$
is soluble  and     $ {\cal L}_{\mathfrak{F}}(G)={\cal
L}_{\Delta}(G)$, then $D$
is
 an abelian  Hall subgroup of odd order of $G$, every element of $G$ 
induces a power automorphism in  
 $D/\Phi (D)$ and every chief factor of $G$ below $D$ is cyclic. }

(ii) {\sl Let    $G$ be  soluble and   $\Delta$    the set of all 
central  chef factors $H/K$  of $G$, that is,  $H/K\leq Z(G/K)$.    If ${\cal
L}_{\mathfrak{N}}(G) ={\cal L}_{\Delta}(G)$, then 
 every element of $G$ induces a power automorphism in $G^{\mathfrak{N}}$. }

Now we consider some applications of Theorem 1.4 in the theory of generalized
$T$-groups. Firstly recall that $G$ is said to be a   \emph{$T$-group}
(respectively
    \emph{$PT$-group},   \emph{$PST$-group}) if  every subnormal subgroup of
$G$ is normal  (respectively permutable, $S$-permutable) in $G$.

Theorem 1.4  allows us to give a new characterization of soluble  $PST$-groups.

{\bf Theorem    1.5.}    {\sl Suppose that  $G$ is soluble. Then $G$ is a
$PST$-group  if and only if      ${\cal
L}_{\mathfrak{N}}(G) ={\cal L}_{\Delta}(G)$, where $\Delta$  is  the set of all 
central  chef factors  of $G$.       }

Since clearly ${\cal L}_{\mathfrak{N}}(G)\subseteq {\cal L}_{sn}(G)$ and, 
in the general case, 
 the lattices ${\cal L}_{\mathfrak{N}}(G)$  and ${\cal L}_{sn}(G)$ do not coincide,
    Theorem 1.5  allows us to strengthen  the following known result.

{\bf Corollary 1.6}  (Ballester-Bolinches,  Esteban-Romero \cite {aust}).
 {\sl If   $G$  is soluble and for every subnormal subgroup $A$  of $G$ we have
 $A/A_{G}\leq Z_{\infty}(G/A_{G})$, then  $G$ is a
$PST$-group.}

 From Theorem 1.4 we get  also the following
well-known result.

{\bf Corollary    1.7 } (Zacher  \cite[ Theorem 2.1.11]{prod}). {\sl    If
 $G$ is soluble  $PT$-group, then
$G$ has
an abelian normal Hall subgroup $D$ of odd order such that $G/D$ is
nilpotent  and every element of $G$ induces a power automorphism in $D$. }

From Theorem 1.5 and Theorem 2.1.11 in \cite{prod} we get   the
following

{\bf Corollary    1.8.}    {\sl Suppose that  $G$ is soluble. Then $G$ is a
$PT$-group  if and only if     ${\cal
L}_{\mathfrak{N}}(G)={\cal L}_{\Delta}(G)$, where  $\Delta$  is  the set of all 
central  chef factors $H/K$  of $G$, and  every
 two  subgroups $A$ and $B$  of any
  Sylow subgroup of $G$ are permutable, that is, $AB=BA$.       }

\section{Proof of Theorem 1.1}

 Direct verification shows that the following two lemmas are true.      

{\bf Lemma 2.1.} {\sl Let $N, M$ and $K < H\leq G$ be normal subgroups of $G$, where $H/K$
is a chief factor of $G$. }

(1) {\sl If $N\leq K$, then   $$(H/K)\rtimes (G/C_{G}(H/K))\simeq ((H/N)/(K/N))\rtimes
 ((G/N)/C_{G/N}((H/N)/(K/N)).$$}

(2) {\sl If $T/L$ is a chief factor of $G$ and $H/K$ and $T/L$  are
$G$-isomorphic, then $C_{G}(H/K)=C_{G}(T/L)$ and $$(H/K)\rtimes (G/C_{G}(H/K))\simeq
(T/L)\rtimes (G/C_{G}(T/L)).$$}

{\bf Lemma 2.2.} {\sl Let $\Delta$ be a $G$-closed set of chief factors of $G$. Let 
   $K\leq H$, $K\leq V$, $W\leq V$ and
 $N\leq H$ be normal subgroups of $G$, where  $H/K\leq 
Z_{\Delta}(G)$.    }

(1) {\sl  $KN/K\leq Z_{\Delta}(G)$   if and only if 
 $N/(K\cap N)\leq Z_{\Delta}(G)$.}

(2) {\sl If  $H/N\leq Z _{\Delta}(G)$, 
 then $H/(K\cap N)\leq Z_{\Delta}(G)$. }

(3) {\sl If $V/K\leq   Z _{\Delta}(G)$,  then   $HV/K\leq   
Z_{\Delta}(G)$. }

{\bf Proof of Theorem  1.1.}  Let $A$ and $B$ be subgroups
 of $G$ such that $A, B\in {\cal L}_{\Delta}(G)$ (respectively
  $A, B\in {\cal L}_{\mathfrak{F}}(G)$).

(1) {\sl $A\cap  B\in {\cal L}_{\Delta}(G)$ (respectively
 $A\cap  B\in {\cal L}_{\mathfrak{F}}(G)$).}

First  note that  $(A\cap B)_{G}=A_{G}\cap B_{G}$.
 On the other hand, 
 from
the $G$-isomorphism $$(A^{G}\cap B^{G})/(A_{G}\cap  B^{G})
=(A^{G}\cap B^{G})/(A_{G}\cap B^{G}\cap A^{G})
\simeq A_{G}(B^{G}\cap A^{G})/A_{G}\leq A^{G}/A_{G}$$  
 we get that 
$(A^{G}\cap B^{G})/(A_{G}\cap B^{G})\leq Z_{\Delta}(G)$  
  (respectively we get that
 $(A^{G}\cap B^{G})/(A_{G}\cap B^{G}) \in \mathfrak{F}$ since $\mathfrak{F}$ is
normally hereditary).  
  Similarly, we  get that
  $(B^{G}\cap A^{G})/(B_{G}\cap A^{G})\leq Z_{\Delta}(G)$  
  (respectively we get that  $(A^{G}\cap B^{G})/(A_{G}\cap B^{G}) \in \mathfrak{F}$.
   But then we get that   $$(A^{G}\cap B^{G})/((A_{G}\cap  B^{G})\cap 
(B_{G}\cap  A^{G}))=  (A^{G}\cap B^{G})/(A_{G}\cap  B_{G})
 \leq  Z_{\Delta}(G) $$ 
by Lemma 2.2(2)
 (respectively we get that $ (A^{G}\cap B^{G})/(A_{G}\cap B_{G})\in
\mathfrak{F}$  since $\mathfrak{F}$ is  a formation).   
But  $(A\cap
B)^{G}\leq A^{G}\cap B^{G}$, so  $$ (A\cap B)^{G}/(A_{G}\cap 
B_{G})=(A\cap B)^{G}/(A\cap B)_{G} \leq Z_{\Delta}(G)$$
  (respectively we get that
 $(A\cap B)^{G}/(A\cap B)_{G} \in \mathfrak{F}$). Therefore
 $A\cap  B\in {\cal L}_{\Delta}(G)$ (respectively
 $A\cap  B\in {\cal L}_{\mathfrak{F}}(G)$).

(2) {\sl  Statement (ii) holds for $G$}.

The set ${\cal L}_{\mathfrak{F}}(G)$ is partially ordered with respect to
set inclusion and $G$ is the greatest element of ${\cal
L}_{\mathfrak{F}}(G)$. Moreover,
 Claim (1) implies  that for any set 
$\{A_{1}, \ldots , A_{n}\} \subseteq  {\cal L}_{\mathfrak{F}}(G)$ we have that
  $A_{1}\cap  \cdots \cap
 A_{n} \in {\cal L}_{\mathfrak{F}}(G)$. Therefore the set  ${\cal
L}_{\mathfrak{F}}(G)$  is a lattice (a meet-sublattice of ${\cal L}(G)$  \cite[p. 7]{Schm}).

(3) {\sl  Statements (i) and (iii) hold for $G$.}

In view of Claim (1), we need only to show that   $\langle A, B \rangle\in  {\cal
L}_{\Delta}(G)$ (respectively $\langle A, B \rangle\in
 {\cal L}_{\mathfrak{F}}(G)$).
From the $G$-isomorphisms $$A^{G}(A_{G}B_{G})/A_{G}B_{G}
\simeq A^{G}/(A^{G}\cap A_{G}B_{G})=A^{G}/ A_{G}(A^{G}\cap
 B_{G})\simeq (A^{G}/A_{G})/(A_{G}(A^{G}\cap B_{G})/A_{G})$$
 we get that  $A^{G}(A_{G}B_{G})/A_{G}B_{G}\leq Z_{\Delta}(G)$ 
 (respectively we get that $A^{G}(A_{G}B_{G})/A_{G}B_{G}\in \mathfrak{F}$ since
 $\mathfrak{F}$ is closed under taking homomorphic images).
  Similarly, we can get that
  $B^{G}(A_{G}B_{G})/A_{G}B_{G}\leq Z_{\Delta}(G)$ 
      (respectively we get that $B^{G}(A_{G}B_{G})/A_{G}B_{G}\in
\mathfrak{F}.)$
  Moreover,  $$A^{G}B^{G}/A_{G}B_{G}
=(A^{G}(A_{G}B_{G})/A_{G}B_{G})(B^{G}(A_{G}B_{G})/A_{G}B_{G})$$ and so
$A^{G}B^{G}/A_{G}B_{G}\leq Z_{\Delta}(G)$ 
 by Lemma 2.2(3) (respectively we have $A^{G}B^{G}/A_{G}B_{G} \in  \mathfrak{F}$ since
 $\mathfrak{F}$ is a   Fitting formation).

Next   note that  $\langle A, B \rangle ^{G}= A^{G}B^{G}$  and
  $A_{G}B_{G}\leq \langle A, B \rangle _{G}$.  Therefore  $ {\langle A, B 
\rangle} ^{G}/{\langle A, B \rangle} _{G}   \leq    
 Z_{\Delta}(G)$ (respectively we get that
$ {\langle A, B \rangle} ^{G}/{\langle A, B \rangle} _{G}   \in \mathfrak{F}$ since
$\mathfrak{F}$ is closed under taking homomorphic images). Hence
 $\langle A, B \rangle\in  {\cal
L}_{\Delta}(G)$ (respectively $\langle A, B \rangle\in
 {\cal L}_{\mathfrak{F}}(G)$).     
The theorem is proved.

\section{Proofs of Theorems 1.4 and 1.5}

\

 {\bf Remark 3.1.}   If $G\in \mathfrak{F}$, where  $\mathfrak{F}$ is a 
formation, then every chief factor of $G$ is $\mathfrak{F}$-central in $G$ 
by  well-known Barnes-Kegel's result \cite[IV,  1.5]{DH}. On 
the other hand,  if   $\mathfrak{F}$ is a saturated formation  and every 
chief factor of $G$ is  $\mathfrak{F}$-central in $G$, then  
$G\in \mathfrak{F}$  by \cite[17.14]{14}.

{\bf Proof of Theorem 1.4.}    (i) Assume that this  is false and let
 $G$ be a counterexample of minimal order. 
 Let $D=G^{\mathfrak{F}}$ be the $\mathfrak{F}$-residual
 of $G$   and  let  $R$ be  a minimal normal subgroup of
$G$.

(1) {\sl The Statement (i)  holds for $G/R$}.

Let $\Delta ^{*}$ be the set of all $ \mathfrak{F}$-central chief factors of $G/R$. 
 By  Proposition
2.2.8   in   \cite{15},  $(G/R)^{\mathfrak{F}}
=RG^{\mathfrak{F}}/R= RD/R  \simeq D/(D\cap R)$
 is soluble. Now let $A/R\in  {\cal
L}_{\mathfrak{F}}(G/R)$. Then 
 from the $G$-isomorphism $$A^{G}/A_{G}\simeq
 (A^{G}/N)/(A_{G}/N)=(A/N)^{G/N}/(A/N)_{G/N}
$$  we get that $A^{G}/A_{G}\in 
\mathfrak{F}$, so $A\in    {\cal
L}_{\mathfrak{F}}(G)$ and hence we have $A\in {\cal L}_{\Delta}(G)$ by 
hypothesis, that is, 
$A^{G}/A_{G}\leq Z_{\Delta}(G)$. It follows that $$ (A/N)^{G/N}/(A/N)_{G/N}
\leq Z_{\Delta ^{*}}(G/R)$$  by Lemma 2.1(1).   Hence $A/R\in  {\cal
L}_{\Delta ^{*}}(G/R)$.   
  Therefore the hypothesis holds
for $G/R$, so we have (1) by the choice of $G$.

(2) {\sl $D$ is nilpotent.}

Assume that this is false. 
Claim (1) implies that $(G/R)^{\mathfrak{F}}=RD/R\simeq D/(R\cap D)$
is nilpotent. Therefore if $R\nleq D$, then $D\simeq D/(R\cap D)=D/1$ is nilpotent.
Consequently    every minimal normal subgroup $N$  of
$G$ is contained in $D$ and $D/N$ is nilpotent. Hence $R$ is abelian. If $N\ne R$, then
 $D\simeq D/1=D/((R\cap D)\cap (N\cap D))$ is
nilpotent. Therefore  $R$ is the unique minimal normal subgroup
of $G$  and  $R\nleq \Phi (G)$ 
 by \cite[Ch. A,
13.2]{DH}. Hence  $R=C_{G}(R)$ by   \cite[Ch. A, 15.6]{DH}.  If $|R|$ is  a prime, then
 $G/R=G/C_{G}(R)$ is cyclic and so $R=D$ is nilpotent.
Thus   $|R|$ is  not a  prime. Let $V$ be a maximal subgroup of $R$.
 Then $V_{G}=1$ and  $V^{G}=R\in {\cal
L}_{\mathfrak{F}}(G)$ since $\mathfrak{F}$ contains all nilpotent groups.
Therefore  $V\in {\cal
L}_{\Delta}(G)$. Hence $V^{G}/V_{G}=R/1$ is $\mathfrak{F}$-central in $G$ and so
 $G/R=G/C_{G}(R)=G/D$, which implies that $D=R$ is nilpotent, a contradiction.
 Therefore we have (2).

(3) {\sl Every subgroup $V$ of $D$ containing $\Phi (D)$
 is normal in $G$.}

Let $V/\Phi (D)$ be a  maximal subgroup of $D/\Phi (D)$. Suppose that $V/\Phi (D)$ is not
normal in $G/\Phi (D)$. Then $V^{G}=D$ and $V\in  {\cal
L}_{\mathfrak{F}}(G)= {\cal L}_{\Delta}(G)$ by Claim (2). Hence
 $D/V_{G}\leq Z_{\Delta}(G)$ and so $G/V_{G}\in
\mathfrak{F}$ by Remark 3.1. But then $D\leq V_{G} < D$.  This
contradiction shows that $V/\Phi (D)$ is normal in $G/\Phi (D)$.
Since
$D/\Phi (D)$ is the direct product of
elementary abelian Sylow subgroups of $D/\Phi (D)$, every subgroup of $D/\Phi (D)$ can be
written as the intersection of some maximal subgroups of $D/\Phi (D)$.
Hence we have (3).

 (4)   {\sl Every   chief factor of $G$ below $D$ is cyclic}
 (This follows from Claim (3) and  Theorem 6.7 in \cite[IV]{DH}).

(5) {\sl  $D$ is a Hall subgroup of $G$. }

 Suppose
that this is false and let $P$ be a   Sylow  $p$-subgroup of $D$ such
that $1 < P < G_{p}$ for some prime $p$ and some Sylow $p$-subgroup $G_{p}$
of $G$. Then $p$ divides  $|G:D|$.

(a)  {\sl $D=P$ is  a minimal normal subgroup of $G$.}

Let $N$ be a minimal normal subgroup of $G$ contained in $D$. Then $N$ is
a $q$-group for some prime $q$ and $NP/N$ is a Sylow $p$-subgroup of $D/N$.  Moreover,
$D/N=(G/N)^{\mathfrak{F}}$  is a Hall subgroup of $G/N$ by
Claim (1) and $p$ divides $|(G/N):(D/N)|=|G:D|$.  Hence   $N=P$ is a Sylow $p$-subgroup of
$D$. Since   $D$ is nilpotent by Claim  (2),  a $p$-complement $V$  of $D$
is characteristic in $D$ and so it is normal in  $G$.
 Therefore $V=1$ and $D=N=P$.

(b) {\sl If  $R\ne D$, then $G_{p}=D\times R$.
  Hence $O_{p'}(G)=1$ and $R/1$ is $\mathfrak{F}$-central in $G$. }

Indeed, $DR/R\simeq D$ is a Sylow
subgroup of $G/R$ by Claims (1) and (a) and hence   $G_{p}R/R=DR/R$, which
implies that $G_{p}=D(G_{p}\cap R)$. But then we have $G_{p}=D\times R$
since $D < G_{p}$ by Claim (a).     Thus  $O_{p'}(G)=1$.   Finally, from 
the $G$-isomorphism  $DR/D\simeq R$   we get that  
$R/1$ is $\mathfrak{F}$-central in $G$.

(c) {\sl $D=R\nleq \Phi (G)$ is the unique minimal normal subgroup of $G$.}

Suppose that  $R\ne D$.
Then   $G_{p}=D\times R$ is elementary abelian $p$-group by Claims (a) and (b).
 Hence  $R=\langle a_{1}
\rangle \times \cdots \times \langle a_{t}\rangle  $ for some elements
$a_{1}, \ldots , a_{t}$ of order $p$. On the other hand,
by Claim  (3), $D=\langle a \rangle $, where $|a|=p$. Now let
 $Z=\langle aa_{1} \cdots a_{t}\rangle$. Then $|Z|=p$ and $ZR=DR=G_{p}$ since
 $Z\cap D=1=Z\cap R$ and $|G_{p}:R|=p$.
If 
$Z=Z_{G}$ is normal in $G$, then  from the $G$-isomorphisms
 $DZ/D\simeq Z$ and $DR/D\simeq R$ we  get that  $Z/1$ is $\mathfrak{F}$-central in 
$G$  since  $R/1$ is $\mathfrak{F}$-central in 
$G$ by Claim (b). Hence $G_{p}=ZR\leq   Z_{\Delta}(G)$ by Lemma 2.2(3).
In the case when $Z_{G}=1$ we have 
$Z < Z^{G}\leq Z_{\Delta}(G)$ by hypothesis and so again we get that 
 $G_{p}=ZR\leq Z_{\Delta}(G)$. But then  $G\in
\mathfrak{F}$ by Remark 3.2. This contradiction shows that we have (c).

(d) {\sl $G$ is supersoluble, so $G_{p}$ is normal in $G$}.

Since
$\mathfrak{F}$ is a saturated formation, $D\nleq \Phi (G)$ and so $D=C_{G}(D)$ by Claim (c) and
\cite[Ch. A, 15.6]{DH}. On the other hand, $|D|=p$ by Claims  (4) and (a), so
 $G/D=G/C_{G}(D)$ is cyclic. Hence $G$ is supersoluble and so for some prime $q$ dividing $|G|$
a Sylow $q$-subgroup $Q$ of $G$  is normal in $G$. Claim (b) implies that, in fact, $Q=G_{p}$.
Hence we have (d).

{\sl The final contradiction for Claim (5).}    Since $\Phi
(G_{p})$ is characteristic in $G_{p}$, Claim (d) implies that  $\Phi
(G_{p})$ is normal in $G$ and so $\Phi (G_{p})\leq \Phi (G) =1$. Hence 
 $G_{p}$ is an elementary
 abelian $p$-group and hence $G_{p}=N_{1}\times \cdots \times N_{n}$ for some
 minimal normal subgroups $N_{1},  \cdots ,  N_{n}$ of $G$ by Mascke's Theorem.
 But then $G_{p}=D$ by Claim (c). This contradiction
  completes the proof of (5).

(6) {\sl Every subgroup $H$ of $D$ is normal in $D$.}

If $H_{G}\ne 1$, then $H/H_{G}$ is normal in 
$D/H_{G}=G^{\mathfrak{F}}/H_{G}$ by Claim (1) 
 and  so
$H$ is normal in $D$.  Now suppose that $H_{G}=1$. Then $H^{G}\leq Z_{\Delta}(G)$
by hypothesis and hence $G/C_{G}(H^{G})\in \mathfrak{F}$ by Theorem 17.14 in 
  \cite{14} and Theorem 6.10 in  \cite[IV]{DH}. It follows that
$D\leq  C_{G}(H^{G})$, which implies that $H$ is normal in $D$.

(7)  {\sl $|D|$ is odd.}

Suppose that $2$ divides $|D|$. Then $G$ has a chief
factor  $ D/K$  with $|D/K|=2$ by Claims (2) and (4).
 But then $D/K\leq Z(G/K)$ and  so $G/K\in \mathfrak{F}$ by
 Remark 3.1, which
 implies that  $D\leq K  < D$. This contradiction shows that we have (7).

(8) {\sl The group D is abelian}

In  view of Claims (6) and (7), $D$
 is a Dedekind group of odd order and so  (8) holds.

 From Claims  (3)--(8) it follows that  Statement (i) holds for
$G$, contrary to the choice of $G$. This final contradiction completes 
the proof of (i).

(ii) We have to  show that if  $H$ is
 any subgroup of  $D=G^{\mathfrak{N}}$, then $x\in N_{G}(H)$ for each $x\in G$.
  It is
enough to consider the case when $H$ is a $p$-group for some prime $p$. 
Moreover, in view of Part (i), we can assume that $x$ is a   $p'$-element 
of $G$.

If $H_{G}\ne 1$, then $H/H_{G}\leq  D/H_{G}=(G/H_{G})^{\mathfrak{N}}$ and 
so the hypothesis holds for $(G/H_{G}, H/H_{G})$ (see the proof of Claim (1)). Thus
  $H/H_{G}$ is normal in $G/H_{G}$ by induction, which implies
that $H$ is normal in $G$. Assume that $H_{G}=1$ and so $H^{G}\leq
Z_{\infty}(G)\cap O_{p}(G)$ since $H$ is subnormal in $G$.  But then  $[H, x]=1$.  
The theorem  is  proved.

{\bf Proof of Theorem 1.5.}    First observe   that if  ${\cal
L}_{\mathfrak{N}}(G)={\cal L}_{\Delta}(G)$, then $G$ is a $PST$-group by
 Theorem 1.4 and Theorem 2.1.8 in \cite{prod}.

 Now assume that   $G$ is a soluble $ PST$-group and let $A\in {\cal L}_{\mathfrak{N}}(G)$, that is,
 $A^{G}/A_{G}$ is nilpotent. Then $A$ is subnormal in $G$ and so 
$A/A_{G}\leq Z_{\infty}(G/A_{G})$ by  Corollary
 2 in \cite{aust} (see also Theorem 2.4.4 in \cite{prod}), which implies that 
 $A^{G}/A_{G}\leq 
Z_{\infty}(G/A_{G})$. Hence $A\in  {\cal L}_{\Delta }(G)$, so 
  $ {\cal L}_{\mathfrak{N}}(G)\subseteq {\cal L}_{\Delta}(G)$
The inverse inclusion follows from the fact that 
if $A\in {\cal L}_{\Delta }(G)$, then 
  $A^{G}/A_{G}\leq Z_{\infty} (G/A_{G})\leq F(G/A_{G})$. 
The theorem  is  proved.

\end{document}